\documentclass{article}

\newenvironment{proof}{\begin{trivlist} \item[]
{\bf Proof.}}{\nolinebreak
\hfill \rule{2mm}{2mm} \end{trivlist}}

\newcounter{ctr}

\newcounter{ctr1}

\newcounter{ctr2}

\newenvironment{example}{\begin{list}{}%
{\addtocounter{definition}{1} \item[] {\bf Example \thedefinition}
\setlength{\itemsep}{0ex} \setlength{\topsep}{0ex}
\setlength{\leftmargin}{3em}} \item }{\end{list}}

\newcounter{ctr3}

\newtheorem{definition}{Definition}[section]    
\newtheorem{theorem}[definition]{Theorem}
\newtheorem{lemma}[definition]{Lemma}
\newtheorem{corollary}[definition]{Corollary}
\newtheorem{proposition}[definition]{Proposition}
\newtheorem{remark}[definition]{Remark}

\newcommand{\NN}{\hbox{{\sf I}\kern-.15em\hbox{\sf N}}}
\newcommand{\ZZ}{\hbox{{\sf Z}\kern-.77em\hbox{\sf Z}\kern.15em}}
\newcommand{\LL}{\hbox{{\sf L}\kern-.77em\hbox{\sf L}\kern.15em}}
\newcommand{\QQ}{\hbox{{\sf I}\kern-.4em\hbox{\sf Q}}}
\newcommand{\RR}{\hbox{{\sf I}\kern-.15em\hbox{\sf R}}}
\newcommand{\CC}{\hbox{{\sf I}\kern-.4em\hbox{\sf C}}}
\newcommand{\TT}{\hbox{{\sf I}\kern-.4em\hbox{\sf T}}}
\newcommand{\DD}{\hbox{{\sf I}\kern-.4em\hbox{\sf D}}}

\newcommand{\sgn}{\,{\rm sgn}\, }

\newcommand{\trace}{\, {\rm trace}\, }

\newcommand{\ud}{\, {\rm d} \kern-.025em }

\newcommand{\mod}[1]{\left| \kern.05em #1 \kern.05em \right|}
\newcommand{\norm}[1]{\left\| \kern.05em #1 \kern.05em \right\|}
\newcommand{\inner}[1]{\left\langle \kern.05em #1 \kern.05em \right\rangle }

\newcommand{\pick}[2]{\renewcommand{\arraystretch}{0.6}
\left( \kern-.4em \begin{array}{c} #1 \\ #2 \end{array} \kern-.4em \right) }

\usepackage{amsfonts}
\begin{document}
\begin{center}{\Large Integrable operators and squares of Hankel matrices}

Andrew McCafferty

\emph{Department of Mathematics and Statistics, Lancaster University, Lancaster, LA1 4YF}
(a.mccafferty@lancaster.ac.uk)
\end{center}
\begin{abstract}
In this note, we find sufficient conditions for an operator with kernel
of the form \\$A(x)B(y)-A(x)B(y)/(x-y)$ (which we call a Tracy--Widom
type operator) to be the square of a Hankel operator. We consider two
contexts: infinite matrices on $\ell^2$, and integral operators on the
Hardy space $H^2(\mathbb{T})$. The results can be applied to the
discrete Bessel kernel, which is significant in random matrix theory.
\end{abstract}
\emph{Keywords}: discrete-time Lyapunov equation, Tracy--Widom operator,
Hankel operator

\emph{2000 Mathematics subject classification}: 47B35, 15A52
\section{Introduction}
In random matrix theory it is natural (see, e.g. \cite{tracywidom}) to
consider integrable operators $T$, where the kernel of
$T$ is
\begin{equation}
\sum_{j=1}^{n}\frac{A_j(z)B_j(w)}{z-w},
\end{equation}
and $\sum_{j=1}^{n}A_j(z)B_j(z)=0$. Here we are concerned with a
special class of such operators, namely those with kernel of the form
\begin{equation}
\label{kernel}
K(x,y)=\frac{A(x)B(y)-A(y)B(x)}{x-y}\quad (x\neq y),
\end{equation}
which we shall refer to as \emph{Tracy--Widom} operators.
The variables $x$ and $y$ may be non-negative integers, as in the
discrete kernels considered in section 2, continuous real parameters,
as in e.g. \cite{blower}, or may live on the circle, as in section 3. We look for conditions under which these operators can be expressed as
$\Gamma^2$ or $\Gamma^*\Gamma$, where $\Gamma$ is a Hankel
operator. In particular we recover a result of Borodin \emph{et al} \cite{borodin}, showing that the discrete Bessel kernel can be written
as
\begin{equation}
\sqrt{\theta}\frac{J_x(2\sqrt{\theta})J_{y+1}(2\sqrt{\theta})-J_y(2\sqrt{\theta})J_{x+1}(2\sqrt{\theta})}{x-y}=\sum_{k=0}^{\infty}J_{x+k+1}(2\sqrt{\theta})J_{y+k+1}(2\sqrt{\theta}).
\end{equation}
We can then read off information about $K$ from knowledge of the
Hankel operator $\Gamma$. For example, a trace
formula follows immediately, and the spectrum of $K$ can be calculated from
the spectrum of $\Gamma$ (which in many cases is easier to
calculate). Megretski, Peller and Treil \cite{pmt} have characterised the
self-adjoint bounded linear operators that are unitarily equivalent to
Hankel operators: we apply their results to gain spectral information
about the operators $K$.
\section{Discrete integrable operators}
Define $\mathbb{N}_0=\mathbb{N}\cup\{0\}$. We consider infinite matrices with kernel
$K(x,y)$, where $K(x,y)$ is defined by (\ref{kernel}). Recall that a Hankel matrix $\Gamma_\phi=[\phi(m+n)]_{m,n\geq
0}$ with $(\phi(k))\in\ell^2$  has square
\begin{equation}
\Gamma_\phi^2=\left[\sum_{k=0}^{\infty}\phi(m+k)\phi(n+k)\right]_{m,n=0}^{\infty}.
\end{equation}
Nehari's theorem (see, e.g. \cite[p. 3]{peller}) states that $\Gamma_\phi$ is
a bounded operator on $\ell^2(\mathbb{N}_0)$ if
and only if $(\phi(n))$ are the positive Fourier coefficients of some
function in $L^\infty(\mathbb{T})$. We write the kernel $K(x,y)$ in matricial form,
\begin{eqnarray}
\nonumber
K(x,y)&=&\frac{1}{x-y}\langle F{\bf a}(x),{\bf a}(y)\rangle,\quad
(x\neq y)\\
{\bf
  a}(x)&=&\left[\begin{array}{c}A(x)\\B(x)\end{array}\right],\quad
  F=\left[\begin{array}{cc}0 & -1\\1 & 0\end{array}\right],
\label{k_definition}
\end{eqnarray}
and look for sufficient conditions under which we can construct a function $\phi:\mathbb{N}_0\rightarrow \CC$ with
$(\phi(j))\in l^2$, such that
\begin{equation}
K(x,y)=\sum_{k=0}^{\infty}\phi(x+k)\phi(y+k),\quad (x\neq y).
\end{equation}
\begin{definition}
Let $S$ be the shift operator on $\ell^2(\mathbb{N}_0)$, so that $Sf(x)=f(x-1)$
(where we define $f(-1)=0)$, and let $R$ be the adjoint shift
operator $Rf(x)=f(x+1)$. The forward difference operator
$\Delta$ is defined by
$\Delta f(x)=f(x+1)-f(x)$. Notice that $\Delta=R_x-I$. Where there are
several variables, we write $R_x$, $\Delta_y$ and so on.
\end{definition}
As usual, $A^T$ is the transpose of a matrix $A$, while $B^*$ denotes
the adjoint of an operator $B$.
\begin{lemma} (Lyapunov equation)
Suppose that $R$ and $B$ are bounded linear operators on $\ell^2$ such that
\[
\sum_{j=0}^{\infty}\left\langle
  R^jBB^{*}(R^{*})^j\xi,\xi\right\rangle<\infty\quad\textrm{
  for all } \xi\in \ell^2,
\]
so that the series
\[
K=\sum_{j=0}^{\infty}R^jBB^{*}(R^*)^j
\]
is convergent in the weak operator topology. Then
\begin{equation}
\label{operator_lyap}
K-RKR^{*}=-BB^{*}.
\end{equation}
\end{lemma}
\begin{proof}
Clear from calculation of the left hand side of (\ref{operator_lyap}).
\end{proof}
In the following Lemma, we state explicitly the specialisation of the
above result to discrete kernels.
\begin{lemma}
\label{lemma1}
Let $\Phi(x,y)$ be any function $\Phi:\mathbb{N}_0^2\rightarrow\mathbb{C}$, and suppose
$\phi:\mathbb{N}_0\rightarrow\mathbb{C}$ is such that $(\phi(j))\in\ell^2$. Then
\begin{equation}
\label{hankelsquared}
\Phi(x,y)=\sum_{k=0}^{\infty}\phi(x+k)\phi(y+k)\quad\textrm{for all }x,y\in\mathbb{N}_0
\end{equation}
if and only if
\begin{equation}
\label{lyap_eq}
(\Delta_xS_y+\Delta_y)\Phi(x,y)=-\phi(x)\phi(y)\quad\textrm{for all }x,y\in\mathbb{N}_0
\end{equation}
and
\begin{equation}
\label{vanishing}
\Phi(x,y)\rightarrow 0 \quad \textrm{as}\quad x\textrm{ or } y\rightarrow\infty.
\end{equation}
\end{lemma}
\begin{proof}
Suppose (\ref{hankelsquared}) holds. Then we have
\begin{eqnarray}
(\Delta_xS_y+\Delta_y)\sum_{k=0}^{\infty}\phi(x+k)\phi(y+k)&=&(S_xS_y-I)\sum_{k=0}^{\infty}\phi(x+k)\phi(y+k)\nonumber\\
&=&\sum_{k=0}^{\infty}\left(\phi(x+k+1)\phi(y+k+1)-\phi(x+k)\phi(y+k)\right)\nonumber\\
&=&-\phi(x)\phi(y),
\end{eqnarray}
so that (\ref{lyap_eq}) holds. By the Cauchy-Schwarz inequality, and since $(\phi(j))\in\ell^2$, we have
\begin{eqnarray}
\Phi(x,y)&=&\sum_{k=0}^{\infty}\phi(x+k)\phi(y+k)\nonumber\\
&\leq& \left(\sum_{k=0}^{\infty}\phi(x+k)^2\right)^{1/2}\left(\sum_{k=0}^{\infty}\phi(y+k)^2\right)^{1/2}\rightarrow
0\quad\textrm{as}\quad x\textrm{ or }y\rightarrow\infty,\label{cauchy_schwarz}
\end{eqnarray}
which is condition (\ref{vanishing}). Conversely, suppose that we have (\ref{lyap_eq}) and
(\ref{vanishing}),  and let
\begin{equation}
G(x,y)=\Phi(x,y)-\sum_{k=0}^{\infty}\phi(x+k)\phi(y+k).
\end{equation}
By (\ref{lyap_eq}), we have
\[
(\Delta_xS_y+\Delta_y)G(x,y)=0\quad\textrm{for all }x,y\in\mathbb{N}_0,
\]
so that $G(x,y)=G(x+1,y+1)$ for all $x,y\in\mathbb{N}_0$. We then use
the hypothesis (\ref{vanishing}) and the estimate in
(\ref{cauchy_schwarz}) to show that $G(x,y)\rightarrow 0$ as $x$ or
$y\rightarrow\infty$, and hence that
$G$ is identically zero for all non-negative integers $x$ and $y$, so
that (\ref{hankelsquared}) holds.
\end{proof}
\begin{theorem}
\label{hankel_prop}
Let $K(x,y)$ be as defined in (\ref{k_definition}), with $({\bf
  a}(x))_{x=0}^{\infty}=([A(x),B(x)]^T)_{x=0}^{\infty}$ a
sequence of $2\times 1$ real vectors such that
\begin{equation}
\sum_{x\geq
  0}\norm{{\bf a}(x)}^2<\infty.
\end{equation}
Suppose that there exists a
sequence of $2\times 2$ real matrices $S_x$ such that ${\bf a}(x+1)=S_x{\bf a}(x)$
for all $x\in\mathbb{N}_0$ and that
\begin{equation}
C=\frac{S_y^{T}FS_x-F}{x-y}
\end{equation}
is a constant matrix. Then $C$ is symmetric. Suppose further that $C$ has
eigenvalues $\lambda\in\RR\setminus\{0\}$ and $0$, and let
$[\alpha,\beta]^T$ be a real unit eigenvector corresponding to $\lambda$. Then
\begin{equation}
\label{k_as_hankel}
K(x,y)=-\sgn(\lambda)\sum_{k=0}^{\infty}\phi(x+k)\phi(y+k)
\qquad \textrm{for} \quad x,y\in\mathbb{N}_0\quad (x\neq y),
\end{equation}
where \begin{equation}
\phi(x)=\mod{\lambda}^{1/2}\left(\alpha
    A(x)+\beta B(x)\right)
\end{equation} 
and $(\phi(x))\in\ell^2$.
\end{theorem}
\begin{proof}
We set
\begin{equation}
C=\frac{S_y^{T}FS_x-F}{x-y}\quad
\end{equation}
where $C$ is constant by hypothesis, so that we can exchange the roles of $x$ and
$y$, and find that $C^T=C$. We have, for $x\neq y$,
\begin{eqnarray}
(\Delta_x
S_y+\Delta_y)K(x,y)&=&(S_xS_y-I)\frac{1}{x-y}\left\langle F{\bf a}(x),{\bf a}(y)\right\rangle\nonumber\\
&=&S_x\frac{1}{x-y-1}\left\langle F{\bf a}(x),S_y{\bf a}(y)\right\rangle-\frac{1}{x-y}\left\langle F{\bf a}(x),{\bf a}(y)\right\rangle\nonumber\\
&=&\frac{1}{x-y}\left\langle FS_x{\bf a}(x),S_y{\bf a}(y)\right\rangle-\frac{1}{x-y}\left\langle F{\bf a}(x),{\bf a}(y)\right\rangle\nonumber\\
&=&\frac{1}{x-y}\left\langle (S_y^TFS_x-F){\bf a}(x),{\bf a}(y)\right\rangle\nonumber\\
&=&\left\langle C{\bf a}(x),{\bf a}(y)\right\rangle.
\end{eqnarray}
 Since $C$ is real and symmetric, and by hypothesis has eigenvalues
 $\lambda\neq 0$ and $0$, there exists a real orthogonal matrix $U$ of unit eigenvectors such that
\begin{equation}
U^TCU=\left[\begin{array}{cc} \lambda & 0\\ 0 & 0\end{array}\right].
\end{equation}
We have
\begin{eqnarray}
(\Delta_xS_y+\Delta_y)K(x,y)&=&\langle C{\bf a}(x),{\bf a}(y)\rangle\nonumber\\
&=&\left\langle U\left[\begin{array}{cc} \lambda & 0\\ 0 & 0\end{array}\right]U^T{\bf a}(x),{\bf a}(y)\right\rangle\nonumber\\
&=&\lambda\left\langle\left[\begin{array}{cc} 1 & 0\\ 0 & 0\end{array}\right]U^T{\bf a}(x),U^T{\bf a}(y)\right\rangle\nonumber\\
&=&\lambda\left\langle\left[\begin{array}{cc} 1 & 0\\ 0 & 0\end{array}\right]U^T{\bf a}(x),\left[\begin{array}{cc} 1 & 0\\ 0 & 0\end{array}\right]U^T{\bf a}(y)\right\rangle\nonumber\\
&=&\sgn(\lambda)\phi(x)\phi(y),
\end{eqnarray}
where
\begin{equation}
\left[\begin{array}{c}\phi(x)\\0\end{array}\right]=\left[\begin{array}{cc}
    \mod{\lambda}^{1/2} & 0\\ 0 & 0\end{array}\right]U^T{\bf a}(x).
\end{equation}
Note that $(\phi(x))\in\ell^2$ by the condition $\sum_{x\geq
  0}\norm{{\bf a}(x)}^2<\infty$, since $U$ is a constant matrix. It is
  also clear that $K(x,y)\rightarrow 0$ as $x$ or $y\rightarrow\infty$, by the same condition on
  ${\bf a}(x)$. We now let $[\alpha,\beta]^T$ be a real unit eigenvector of $C$ corresponding to $\lambda$, and the result follows by Lemma \ref{lemma1}.
\end{proof}
\begin{corollary}
\label{corollary}
Let $K(x,y)$ be as defined in (\ref{k_definition}), with $({\bf
  a}(x))_{x=0}^{\infty}=([A(x),B(x)]^T)_{x=0}^{\infty}$ a
sequence of $2\times 1$ real vectors such that \begin{equation}
\sum_{x\geq
  0}\norm{{\bf a}(x)}^2<\infty.
\end{equation}
 Suppose that ${\bf
  a}(x+1)=(Lx+M){\bf a}(x)$ (for all $x\in\mathbb{N}_0$), where $L$
  and $M$ are real constant $2\times 2$ matrices that satisfy
\[
\left\{
\begin{array}{l}
\det L=0\\
\det M=1\\
M^TFL \textrm{ is symmetric, and has eigenvalues } \lambda\in\RR\setminus\{0\} \textrm{ and } 0.
\end{array}
\right.
\]
Let $[\alpha,\beta]^T$ be a real unit eigenvector of $M^TFL$ corresponding
to $\lambda$. Then
\begin{equation}
K(x,y)=-\sgn(\lambda)\sum_{k=0}^{\infty}\phi(x+k)\phi(y+k)
\ \textrm{ for all } x,y\in\mathbb{N}_0\quad (x\neq y),
\end{equation}
where $\phi(x)=\mod{\lambda}^{1/2}\left(\alpha A(x)+\beta B(x)\right)$, and $(\phi(x))\in\ell^2$.
\end{corollary}
\begin{proof}
We have $M^TFM=F\det M$ (indeed, this is true for any $2\times 2$ matrix) and hence\\
$M^TFM=F$. Likewise $L^TFL=0$. Setting $S_x=Lx+M$ as
in Theorem \ref{hankel_prop}, we now have
\begin{eqnarray}
\frac{S_y^TFS_x-F}{x-y}&=&\frac{(Ly+M)^TF(Lx+M)-F}{x-y}\nonumber\\
&=&\frac{M^TFLx-(M^TFL)^Ty}{x-y} \qquad \textrm{(since } F^T=-F)\nonumber\\
&=&\frac{M^TFL(x-y)}{x-y} \qquad \textrm{(since } M^TFL \textrm{ is symmetric
by hypothesis)}\nonumber\\
&=&M^TFL.
\end{eqnarray}
Hence $C=(S_y^TFS_x-F)/(x-y)$ is a constant matrix. Thus, together with the
summability criterion on the sequence $({\bf a}(x))$, the hypotheses of
Theorem \ref{hankel_prop} are all satisfied, so we have the result.
\end{proof}
\begin{example}
Let $J_\nu(z)$ be the Bessel functions of the first kind of order
$\nu$, and write $J_x=J_x(2\sqrt{\theta})$, where $\theta$ is a
positive real parameter. The discrete Bessel kernel
\begin{equation}
J(x,y;\theta)=\sqrt{\theta}\frac{J_xJ_{y+1}-J_yJ_{x+1}}{x-y}
\end{equation}
arises in the study of various discrete-variable random matrix
models, as in \cite{johannson} and \cite{borodin}. Note that $J_x$ is
an entire function of $x$, so that $J(x,x;\theta)$ is well-defined via
L'Hopital's rule. In the notation of
Corollary \ref{corollary}, we take
\begin{equation}
{\bf{a}}(x)=\left[\begin{array}{c}\sqrt{\theta}J_x\\J_{x+1}\end{array}\right].
\end{equation}
The
standard formula (see \cite[p. 379]{whit_watson})
\begin{equation}
e^{i2t\sin{\theta}}=J_0(2t)+2\sum_{m=1}^{\infty}J_{2m}(2t)\cos 2m\theta+2i\sum_{m=1}^{\infty}J_{2m-1}(2t)\sin(2m-1)\theta
\end{equation}
and Parseval's identity can be used to show that
$J_0(2t)^2+2\sum_{m=1}^{\infty}J_{m}(2t)^2=1$ for all real $t$, and hence that the
sequence $(J_{x})_{x=0}^{\infty}$ is square summable. Thus the condition $\sum_{x\geq
  0}\norm{{\bf a}(x)}^2<\infty$ is satisfied. The 3-term recurrence relation for the Bessel functions
\begin{equation}
J_{x+2}(2z)-\frac{x+1}{z}J_{x+1}(2z)+J_x(2z)=0
\end{equation}
becomes
\begin{equation}
{\bf
  a}(x+1)=\left[\begin{array}{c}\sqrt{\theta}J_{x+1}\\J_{x+2}\end{array}\right]=\left[\begin{array}{cc}0
  & \sqrt{\theta}\\\frac{-1}{\sqrt{\theta}} & \frac{x+1}{\sqrt{\theta}}\end{array}\right]\left[\begin{array}{c}\sqrt{\theta}J_x\\J_{x+1}\end{array}\right],
\end{equation}
and so we have ${\bf{a}}(x+1)=(Lx+M){\bf{a}}(x)$, where
\[
L=\left[\begin{array}{cc}0 & 0\\0 &
    \frac{1}{\sqrt{\theta}}\end{array}\right]
\]
and
\[
M=\left[\begin{array}{cc}0 & \sqrt{\theta}\\\frac{-1}{\sqrt{\theta}} & \frac{1}{\sqrt{\theta}}\end{array}\right].
\]

It is clear that these matrices satisfy $\det L=0$ and $\det M=1$, and we have
\[
M^TFL=\left[\begin{array}{cc}0 & 0\\0 & -1\end{array}\right],
\]
so we pick the unit eigenvector $[\alpha,\beta]^T=[0,1]^T$. Thus, the
function $\phi(x)$ in Corollary \ref{corollary} is $J_{x+1}(2\sqrt{\theta})$, and we recover a result
of Borodin \emph{et al} in \cite{borodin}
\begin{equation}
J(x,y;\theta)=\sum_{k=0}^{\infty}J_{x+k+1}(2\sqrt{\theta})J_{y+k+1}(2\sqrt{\theta}),
\quad x,y\in\mathbb{N}_0,
\end{equation}
without their use of asymptotic formulae for the Bessel functions.
\end{example}
The preceding results are identities of kernels for $x\neq y$. Evidently, the sum in the
right-hand side of (\ref{k_as_hankel}) makes sense for $x=y$, and hence gives one possible
extension of the left-hand side to the case $x=y$. We use the
extension to define an operator $K$ with matrix given by $K(x,y)$.
\begin{proposition}
\label{traceclass}
Suppose that the vector $[A(x),B(x)]^T$ satisfies the conditions of
Theorem \ref{hankel_prop}, so that
$K(x,y)=\sum_{k=0}^{\infty}\phi(x+k)\phi(y+k)$. Suppose also that
$\sum_{n=0}^{\infty}n\phi(n)^2<\infty$. Then the operator
$K$ represented by the matrix $[K(x,y)]_{x,y=0}^{\infty}$ is trace
class and has trace:
\begin{equation}
\trace K=\sum_{x=0}^{\infty}(x+1)\phi(x)^2.
\end{equation}
\end{proposition}
\begin{proof}
The summability condition on $\phi$ ensures that $\Gamma_\phi$ is
Hilbert-Schmidt, which implies that $K=\Gamma_\phi^2$ is trace-class. We have
\begin{equation}
\trace K=\sum_{x=0}^{\infty}K(x,x)=\sum_{x=0}^{\infty}\sum_{k=0}^{\infty}\phi(x+k)^2
\end{equation}
from which the result follows immediately.
\end{proof}
\begin{definition}
For a compact and self-adjoint operator $W$ on a Hilbert space $H$,
the spectral multiplicity function $\nu_W(\lambda):\RR\rightarrow
\{0,1,\ldots\}\cup\{\infty\}$ is given by
\begin{equation}
\nu_W(\lambda)=\dim\{x\in H: Wx=\lambda x\}\quad (\lambda\in\RR).
\end{equation}
\end{definition}
We now give the consequences of a
result of Peller, Megretski and Treil in \cite{pmt} in the case of discrete integrable operators.
\begin{proposition}
\label{spectral}
Suppose that $\Gamma_\phi$ and $K$ are as in Proposition
\ref{traceclass}. Then $\Gamma_\phi$ and $K$ are compact and self-adjoint, and
\renewcommand{\labelenumi}{\roman{enumi}}
  \newcounter{Rcount}
  \begin{list}{(\roman{Rcount})}
    {\usecounter{Rcount}
    \setlength{\rightmargin}{\leftmargin}}
\item $\nu_{K}(0)=0$ or $\nu_K(0)=\infty$;
\item for $\lambda>0$, $\nu_K(\lambda)<\infty$ and
  $\nu_K(\lambda)=\nu_{\Gamma\phi}(\sqrt{\lambda})+\nu_{\Gamma_\phi}(-\sqrt{\lambda})$;
\item if $\nu_K(\lambda)$ is even, then
  $\nu_{\Gamma_\phi}(\sqrt{\lambda})=\nu_{\Gamma\phi}(-\sqrt{\lambda})$;
\item if $\nu_K(\lambda)$ is odd, then $\mod{\nu_{\Gamma\phi}(\sqrt{\lambda})-\nu_{\Gamma_\phi}(-\sqrt{\lambda})}=1$.
\end{list}
\end{proposition}
\begin{proof}
\emph{(i)} follows from Beurling's theorem (see \cite{peller}, page 15), while \emph{(ii)} is
elementary. Peller, Megretski and Treil show in \cite{pmt} that for
any compact and self-adjoint Hankel operator $\Gamma$, the spectral
multiplicity function satisfies
$\mod{\nu_\Gamma(\lambda)-\nu_\Gamma(-\lambda)}\leq 1$. Using this, and
\emph{(ii)}, statements
\emph{(iii)} and \emph{(iv)} follow immediately.
\end{proof}
\begin{remark}
The Carleman operator $\Gamma: L^2(0,\infty)\rightarrow L^2(0,\infty)$ is
given by
\begin{equation}
\Gamma f(x)=\int_{0}^{\infty}\frac{1}{x+t}f(t)\ud t,
\end{equation}
so $\Gamma^2$ has kernel of Tracy-Widom type 
\begin{equation}
\Gamma^2f(u)=\int_{0}^{\infty}\frac{\log u-\log t}{u-t}f(t)\ud
t.
\end{equation}
Carleman showed that $\Gamma$ is a positive self-adjoint Hankel
operator with continuous spectrum $[0,\pi]$ of multiplicity two (see
\cite[p. 442]{peller}), so the Tracy--Widom type operator $\Gamma^2$
has spectrum $[0,\pi^2]$, also of multiplicity two. This contrasts with
\emph{(iii)} and \emph{(iv)} of Proposition \ref{spectral}.
\end{remark}
\section{Integrable operators on $H^2$}
Let $H^2$ be the usual Hardy space on the unit circle $\mathbb{T}$, with
orthonormal basis $\{1,z,z^2,\ldots\}$, and let $R_+:
L^2\rightarrow H^2$ and $R_-:L^2\rightarrow L^2\ominus H^2$ be the
Riesz orthogonal projection operators. We let $M_\phi$ denote
multiplication by $\phi$, and define the Toeplitz operator on $H^2$ with symbol
$\phi$ to be $T_\phi=R_+M_\phi R_+$. Let $J:L^2\rightarrow L^2$ be a flip operator, whose
operation on a function $f\in H^2$ is $Jf(z)=\bar{z}f(\bar{z})$. Note that $J$ maps $H^2$ onto $L^2\ominus
H^2$ (and vice versa) and that $J^2=I$. The
Hankel operator $\Gamma_\phi$ on $H^2$ with symbol $\phi\in
L^{\infty}$ is then
\begin{equation}
\Gamma_\phi=JR_-M_\phi.
\end{equation}
We let the integral operator $W$ on $L^2(\mathbb{T})$ have kernel
\begin{equation}
W(e^{i\theta},e^{i\phi})=\frac{f(e^{i\theta})g(e^{i\phi})-f(e^{i\phi})g(e^{i\theta})}{1-e^{i(\theta-\phi)}},
\end{equation}
where $W$ operates on a function $f\in L^2(\mathbb{T})$ in the usual way:
\begin{equation}
Wf(e^{i\theta})=\frac{1}{2\pi}\int_{\mathbb{T}} W(e^{i\theta},e^{i\phi})f(e^{i\phi})\ud \phi.
\end{equation}
\begin{lemma}
\label{projection_prop}
Suppose that $f,g\in L^{\infty}$ have $\bar{f}=g$. Then $W$ defines a
bounded and self-adjoint operator on $L^2$. Further,
$R_{+}WR_{+}:H^2\rightarrow H^2$ satisfies
\begin{equation}
\label{hankel_equation}
R_{+}WR_{+}=\Gamma_{f}^{*}\Gamma_{f}-\Gamma_{g}^{*}\Gamma_{g}.
\end{equation}
Moreover, when $f$ is continuous, $R_{+}WR_{+}$ is compact.
\end{lemma}
\begin{proof}
The condition $\bar{f}=g$ gives immediately
$W(e^{i\theta},e^{i\phi})=\overline{W(e^{i\phi},e^{i\theta})}$, and so
  $W$ is self-adjoint. It can easily be seen that the Riesz projection
  $R_{+}$ has distributional kernel $1/(1-e^{i(\theta-\phi)})$,
and so $W$ decomposes as
\begin{equation}
W=M_g\left[M_f,R_+\right]-M_f\left[M_g,R_+\right],
\end{equation}
where all the operators are bounded. A simple calculation now shows that
\begin{equation}
R_+WR_+=(T_{gf}-T_gT_f)-(T_{fg}-T_fT_g),
\end{equation}
and we apply the standard formulae
$T_{hk}-T_hT_k=\Gamma_{h(\bar{z})}\Gamma_{k(z)}$ and
$\Gamma_{h}^{*}=\Gamma_{\bar{h}(\bar{z})}$ (see \cite[p. 253]{nikolski})
to get equation (\ref{hankel_equation}). The last statement follows by
Hartman's theorem: the Hankel operators on the right-hand side of
(\ref{hankel_equation}) are compact when $f$ is continuous.
\end{proof}
\begin{remark}
We continue functions $f\in L^2$ to harmonic functions on $\mathbb{D}$ by
means of the Poisson kernel, as in \cite[p. 718]{peller}.
\end{remark}
\begin{proposition}
Suppose $f=\bar{g}\in L^{\infty}$, where $g$ is holomorphic inside
$\mathbb{D}$. Then 
\begin{equation}
R_+WR_+=\Gamma^*_{f}\Gamma_f.
\end{equation}
Further, if $R_+WR_+$ has finite rank,
then $f$ is rational.
\end{proposition}
\begin{proof}
Take $f=\bar{g}$ in Lemma \ref{projection_prop} to obtain the first part
of the result. For the second part, note that
\begin{equation}
\textrm{Range}(R_+WR_+)=\textrm{Ker}(\Gamma_f^*\Gamma_f)^\perp=\textrm{Ker}(\Gamma_f)^\perp=\textrm{Range}(\Gamma_f^*)=\textrm{Range}(\Gamma_{\bar{f}(\bar{z})}),
\end{equation}
and apply Kronecker's theorem: $\Gamma_k$ has finite rank if and only
if $k$ is rational, so $\Gamma_{\bar{f}(\bar{z})}$ has finite rank if
and only if $\bar{f}(\bar{z})$ is rational, which implies that $f$ is rational.
\end{proof}
{\bf Acknowledgements}\\
I would like to thank my supervisor, Gordon Blower for his patience in many discussions. Thanks are also due to Steve Power for helpful suggestions. I am funded by the EPSRC's Doctoral Training Account scheme.


\begin{thebibliography}{9}
\bibitem{tracywidom} C.A Tracy and H. Widom, Level-spacing
  distributions and the Airy kernel, Comm. Math. Phys. 159 (1994), 151-174
\bibitem{blower} G Blower Operators associated with Soft and Hard
  Spectral Edges from Unitary Ensembles, J. Math. Anal. Appl. 2007.
\bibitem{borodin} A. Borodin, A. Okounkov, G. Olshanski, Asymptotics
  of Plancherel measures for symmetric groups, J. Amer. Math. Soc.  13
  (2000) 481-515.
\bibitem{pmt}  A.V. Megretskii, V.V. Peller, and S.R. Treil, The
  inverse spectral problem for self-adjoint Hankel operators, Acta
  Math. 174 (1995),  no. 2, 241-309
\bibitem{peller} V.V. Peller, Hankel Operators and Their
    Applications, Springer, New York, 2003.
\bibitem{johannson} K. Johansson, Discrete Orthogonal Polynomial
  ensembles and the Plancherel measure, Ann. Math. (2) 153 (2001), no. 2, 259-296.
\bibitem{power}S.C. Power, Hankel Operators on Hilbert Space,
  Research Notes in Math., 64, Pitman Adv. Publ. Progr.,
  Boston-London-Melbourne, 1982.
\bibitem{whit_watson}E.T Whittaker and G.N. Watson, A Course of Modern
  Analysis, Cambridge University Press, 4th ed. 1963
\bibitem{nikolski} N.K. Nikolski, Operators, Functions and
    Systems: an easy reading, vol. 2 American Mathematical Society, Providence, R.I., 2002.
\end{thebibliography}
\end{document}